\documentclass[12pt]{article}
\usepackage[english]{babel}
\usepackage{amsfonts}
\usepackage{amsmath}
\usepackage{amssymb}
\usepackage{mathrsfs}
\usepackage{latexsym}
\usepackage{multicol}
\usepackage{multirow}
\usepackage{fancybox} 
\usepackage{color}
\usepackage{hyperref}
\usepackage{graphicx} 
\usepackage{caption}
\usepackage{enumitem}
\usepackage{upgreek}
\usepackage[overload]{empheq}
\usepackage{diagbox}
\usepackage{adjustbox}
\usepackage{caption, multirow, makecell}

\usepackage{caption}
\usepackage{graphicx}
\usepackage{upgreek}
\usepackage[overload]{empheq}
\usepackage{float}
\def\be{\begin{equation}}
\def\ee{\end{equation}}

\def\d'{``}

\def\be{\begin{equation}}
\def\ee{\end{equation}}
\def\bea{\begin{eqnarray}}
\def\eea{\end{eqnarray}}

\def\i'{\textrm{i}}

\def\d'{``}
\begin{document}

\title{Lommel functions, Pad\'e approximants and hypergeometric functions.} 
\author{Federico Zullo\thanks{DICATAM, Universit\`a di Brescia, Brescia, Italy \& INFN, Sezione di Milano-Bicocca, Milano, Italy.}}
\date{}

\maketitle 

\begin{abstract}
We consider the Lommel functions $s_{\mu,\nu}(z)$ for different values of the parameters $(\mu,\nu)$. We show that if $(\mu,\nu)$ are half integers, then it is possible to describe these functions with an explicit combination of polynomials and trigonometric functions. The polynomials turn out to give Pad\'e approximants for the trigonometric functions. Numerical properties of the zeros of the polynomials are discussed. Also, when $\mu$ is an integer, $s_{\mu,\nu}(z)$ can be written as an integral involving an explicit combination of trigonometric functions. A closed formula for $_2F_1\left(\frac{1}{2}+\nu,\frac{1}{2}-\nu;\mu+\frac{1}{2};\sin(\frac{\theta}{2})^2\right)$ with $\mu$ an integer is given.
\end{abstract}


\section{Introduction}\label{ra_sec1}
When $\nu^2 \neq (\mu+2k+1)^2, k=0,1,2,...$, the Lommel functions $s_{\mu,\nu}(z)$ are defined by the converging series \cite{NIST}, \cite{W}:
\begin{equation}
s_{\mu,\nu}(z)=\frac{z^{\mu+1}}{(\mu+1)^2-\nu^2}\left(1-\frac{z^{2}}{(\mu+3)^2-\nu^2}+\frac{z^{4}}{\left((\mu+3)^2-\nu^2\right)\left((\mu+5)^2-\nu^2\right)}+\cdots\right),
\end{equation}
and are particular solutions of the following inhomogeneous Bessel differential equation
\begin{equation}\label{Lom}
z^2\frac{d^2y}{dz^2}+z\frac{dy}{dz}+(z^2-\nu^2)y=z^{\mu+1}.
\end{equation}
Equivalently, $s_{\mu,\nu}(z)$ can be written in terms of hypergeometric $_1F_2$ function as:
\begin{equation}\label{iper12}
s_{\mu,\nu}(z)=\frac{z^{\mu+1}}{(\mu+1)^2-\nu^2}\,_1F_2\left(1;\frac{\mu-\nu+3}{2},\frac{\mu+\nu+3}{2};-\frac{z^2}{4}\right)
\end{equation}
Among the many properties of these function (see e.g. \cite{NIST} and \cite{W} and references therein for a complete list), we remember here the recurrences solved by $s_{\mu,\nu}(z)$ and their derivatives:
\begin{equation}\label{recu}\begin{split}
&s_{\mu+2,\nu}(z)=z^{\mu+1}- \left((\mu+1)^2-\nu^2\right)s_{\mu,\nu}(z),\\
& \frac{ds_{\mu,\nu}}{dz}\pm \frac{\nu}{z}s_{\mu,\nu}(z)=(\mu\pm \nu -1)s_{\mu-1,\nu\mp 1}(z)
\end{split}\end{equation}
The Lommel functions have many applications in mathematical physics and applied sciences: the interested reader can look at \cite{NIST} and \cite{Z1} and references therein.

This paper is the second of a series of paper dedicated to the Lommel functions. In the first paper \cite{Z1} the distribution of the zeros of the functions $s_{\mu,\nu}(z)$ has been analyzed by making use of an integral representation of $s_{\mu,\nu}(z)$ involving the hypergeometric $\,_1F_2$ functions and old results due to P\'olya \cite{P} about the zeros of functions defined by trigonometric integrals. %by showing how these functions possess only simple zeros (apart the zero $z=0$ of order $\mu+1$ when $\mu$ is an integer). 
In this paper we will discuss the relevant special cases whit $(\mu,\nu)=(m+1/2, n+1/2)$ with $(m,n)$ integers and the special cases $\mu=n$ with $n$ an integer. In section 2 we will remember some results from \cite{Z1} that will be useful in the successive part of the work. In section 3 we will show how it is possible to explicitly write an expression involving polynomials and the classical trigonometric functions for the functions $s_{m+\frac{1}{2},n +\frac{1}{2}}(z)$. The coefficients of the polynomials can be explicitly described in terms of the derivatives of the Legendre polynomials $P_{2n}(t)$ and of the associated Legendre $P_{2n+1}^{1}(t)$ evaluated in $t=0$ and $t=1$.  Further, the ratios of the polynomials are Pad\'e approximants of order $z^{m+n+1}$ to the trigonometric functions $\sin(z)$ and $\cos(z)$. In section 4, after giving the explicit form of the coefficients of the polynomials of section 3, we will give some numerics about them. In section 5 we will write $s_{n,\nu}(z)$ as an integral involving only trigonometric functions: the link with the integral of $s_{\nu,\mu}(z)$ involving the  hypergeometric $\,_1F_2$ function will give an explicit formula for the function $_2F_1\left(\frac{1}{2}+\nu,\frac{1}{2}-\nu;n+\frac{1}{2};\sin(\frac{\theta}{2})^2\right)$ in terms of $\theta$.

\section{Integral representations and zeros.}
In \cite{Z} it has been noticed that the function $s_{0,v}(z)$ can be represented by the following integral:
\begin{equation}\label{int}
s_{0,\nu}(z)=\frac{1}{1+\cos(\pi\nu)}\int_0^\pi \sin(z\sin(t))\cos(\nu t)dt.
\end{equation}
The consequences of these type of integral representation are many and noteworthy in the study of the properties of these functions. For example, thanks to a Theorem by P\'olya, it has been shown in \cite{Z} that $s_{0,\nu}(z)$, for $|\nu|<1$, possesses only real zeros: these zeros are simple and the intervals $(k\pi, (k+1)\pi)$, $k=0, 1, 2, \ldots$ contain the non-negative zeros, each interval containing just one zero.  
In \cite{Z1} an integral representation has been given for $s_{\mu,\nu}(z)$ for $\mu >1/2$. In particular it has been shown that the following formula holds:
\begin{equation}\label{intr}
s_{\mu,\nu}(z)=z^{\mu}\int_0^1 \sin(zt)f_{\mu,\nu}(t)dt=z^{\mu}\int_0^1 \frac{\sin(zt)}{(1-t)^{\frac{1}{2}-\mu}}\frac{_2F_1\left(\frac{1}{2}+\nu,\frac{1}{2}-\nu;\mu+\frac{1}{2};\frac{1-t}{2}\right)}{_2F_1\left(\frac{1}{2}+\nu,\frac{1}{2}-\nu;\mu+\frac{1}{2};\frac{1}{2}\right)}dt.
\end{equation} 
The function $f_{\mu,\nu}(t)$, defined by (\ref{intr}), possesses different properties. For example it solves a differential recurrence:
\begin{equation}\label{deriv1}
\frac{df_{\mu,\nu}}{dt}+a_{\mu,\nu}f_{\mu-1,\nu}=0,
\end{equation}
where we set
\begin{equation}\label{amn}
a_{\mu,\nu}\doteq 2\frac{\Gamma\left(\frac{\mu+1+\nu}{2}\right)\Gamma\left(\frac{\mu+1-\nu}{2}\right)}{\Gamma\left(\frac{\mu+\nu}{2}\right.)\Gamma\left(\frac{\mu-\nu}{2}\right)}.
\end{equation}
with $\mu\pm\nu$ are different from an odd negative integer. Another integral representation is the following \cite{Z1}
\begin{equation}\label{smnnn}
s_{\mu,\nu}=\frac{a_{\mu+2,\nu}}{((\mu+1)^2-\nu^2)}z^{\mu+1}\int_0^1 \cos(zt)f_{\mu+1,\nu}(t).
\end{equation}
Again, from these integral representation, and in particular from the monotonicity properties of $f_{\mu,\nu}(t)$, different properties of $s_{\mu,\nu}(z)$ can be derived. For example it has been shown \cite{Z1} that, apart the branch point at $z=0$, the functions $s_{\mu,\nu}(z)$ for $\mu \in (-1/2,1/2) \cap \nu \in (|\mu|, \mu+1)$ possess only real zeros. The zeros are simple and the intervals $(k\pi, (k+1)\pi)$, $k=0, 1, 2, \ldots$ contain the non-negative zeros, each interval containing just one zero. Another result is the following: apart the branch point at $z=0$, the functions $s_{\mu,\nu}(z)$ for $\mu >1/2 \cap \nu \in (0,\mu)$ are positive on the positive real axis, i.e. they possess only complex zeros. Further, it is possible to show that the function $z^{-\mu}\left(a_{\mu,\nu}\cos(\theta)s_{\mu-1,\nu}(z)+\sin(\theta)s_{\mu,\nu}(z)\right)$ for $\mu \in (-1/2,1/2)$ and  $|\nu| \in (|\mu|, \mu+1)$ possesses only real zeros for any $\theta \in [0,\pi]$. The zeros are simple and each of the intervals
\begin{equation}\label{genint}
\left(\left(k-\frac{1}{2}\right)\pi+\theta,\left(k+\frac{1}{2}\right)\pi+\theta\right),\quad k=0,\pm 1, \pm 2, \ldots
\end{equation} 
contain just one zero. In the following, we will make use again of the integral representation (\ref{intr}) to derive further properties of the Lommel functions.

\section{The case of algebraic kernels}
For particular values of the parameters $\mu$ and $\nu$ it is possible to get algebraic functions for the kernel of the integral (\ref{intr}). First of all, let us 
notice that, for $\nu=\frac{1}{2}$, equation (\ref{intr}) gives a well known representation for $s_{\mu,\frac{1}{2}}(z)$ (see e.g. Koumandos and Lamprecht \cite{KL}):
 \begin{equation}\label{intr0}
s_{\mu,\frac{1}{2}}(z)=z^{\mu}\int_0^1 (1-t)^{\mu-\frac{1}{2}}\sin(zt)dt.
\end{equation}
Actually, when $\nu=\frac{2n+1}{2}$, with $n \in \mathbb{N}$, the hypergeometric function collapses in a polynomial, since we get:
\begin{equation}\label{pol}
s_{\mu,n+\frac{1}{2}}(z)=z^{\mu}\int_0^1 \frac{\sin(zt)}{(1-t)^{\frac{1}{2}-\mu}}\frac{_2F_1\left(-n,n+1;\mu+\frac{1}{2};\frac{1-t}{2}\right)}{_2F_1\left(-n,n+1;\mu+\frac{1}{2};\frac{1}{2}\right)}dt.
\end{equation} 
Taking into account the series of the hypergeometric function $_2F_1\left(-n,n+1;\mu+\frac{1}{2};\frac{1-t}{2}\right)$, we can write the following explicit sum for $s_{\mu,\frac{2n+1}{2}}(z)$:
\begin{equation}\label{intex}
s_{\mu,n+\frac{1}{2}}(z)=\frac{\Gamma\left(\frac{2\mu+3+2n}{4}\right)\Gamma\left(\frac{2\mu+1-2n}{4}\right)}{2^{\frac{1}{2}-\mu}n!\sqrt{\pi}}z^\mu\sum_{k=0}^n \frac{(-1)^k(n+k)!}{2^k\Gamma\left(k+\mu+\frac{1}{2}\right)}\binom{n}{k}\int_0^1\sin(zt)\left(1-t\right)^{k+\mu-\frac{1}{2}}.
\end{equation} 
The previous can be also rewritten in terms of products as
\begin{equation}\label{intex1}
s_{\mu,n+\frac{1}{2}}(z)=z^{\mu}\prod_{p=0}^{n-1}\frac{2\mu-4p+2n-1}{2\mu+2p-2n+1}\sum_{k=0}^{n}\prod_{q=0}^{k-1}\frac{(q-n)(q+n+1)}{(q+1)(2\mu+2q+1)}\int_0^1\sin(zt)\left(1-t\right)^{k+\mu-\frac{1}{2}}
\end{equation}

After $s_{\mu,\frac{1}{2}}(z)$, the first two elements of the set of functions $s_{\mu,n+\frac{1}{2}}(z)$, $n\in \mathbb{N}$, are:
\begin{equation}\label{pol1}
s_{\mu,\frac{3}{2}}(z)=z^{\mu}\int_0^1 (1-t)^{\mu-\frac{1}{2}}\left(1+\frac{2}{2\mu-1}t\right)\sin(zt)dt,
\end{equation}
and 
\begin{equation}\label{pol2}
s_{\mu,\frac{5}{2}}(z)=z^{\mu}\int_0^1 (1-t)^{\mu-\frac{1}{2}}\left(1+\frac{6(2\mu-1)}{(2\mu+1)(2\mu-3)}t+\frac{12}{(2\mu+1)(2\mu-3)}t^2\right)\sin(zt)dt.
\end{equation}

From (\ref{intex})  it readily follows that if $\mu=m+\frac{1}{2}$, with $m\in \mathbb{N}$ and $m-n$, $m+n+1$ not negative odd integers, then the integral (\ref{intex}) gives a combination of rational functions of $z$ and trigonometric functions. More precisely we can define
\begin{equation}\label{smn}
s_{\frac{2m+1}{2},\frac{2n+1}{2}}(z)=\frac{\left(A_{m,n}(z)-B_{m,n}(z)\cos(z)-C_{m,n}(z)\sin(z)\right)}{z^{n+\frac{1}{2}}},
\end{equation}
where the three functions $A_{m,n}(z)$, $B_{m,n}(z)$ and $C_{m,n}(z)$ are polynomials in $z$. Indeed, it is not difficult to show by induction that if $R_k(t)$ is a polynomial of degree $k$ then one has
\begin{equation}\label{intsin}
\int \sin(zt)R_k(t)dt=\sum_{j=0}^{\left \lfloor{\frac{k}{2}}\right \rfloor }\frac{\left((-1)^{j+1}z^{k-2j}\cos(zt)R_k^{2j}(t)+(-1)^j z^{k-2j-1}\sin(zt)R_k^{2j+1}(t)\right)}{z^{k+1}},
\end{equation}
where the apex on $R$ stands for the order of the derivative. In order to identify which are the polynomials  $A_{m,n}(z)$, $B_{m,n}(z)$ and $C_{m,n}(z)$ in (\ref{smn}), from (\ref{recu}) we notice that they satisfy the following difference equations with respect to $m$:
\begin{equation}\begin{split}
&A_{m+2,n}(z)+(m+n+2)(m+1-n)A_{m,n}(z)=z^{m+n+2},\\
&B_{m+2,n}(z)+(m+n+2)(m+1-n)B_{m,n}(z)=0,\\
&C_{m+2,n}(z)+(m+n+2)(m+1-n)C_{m,n}(z)=0.
\end{split}\end{equation}
The previous equations show that $C_{m,n}(z)$ and $B_{m,n}(z)$ are proportional, respectively, to $C_{0,n}(z)$ and $C_{1,n}(z)$ and to $B_{0,n}(z)$ and $B_{1,n}(z)$ with coefficients independent of $z$. More precisely one has:
\begin{equation}
\begin{split}
&B_{2m,n}(z)=\frac{(-1)^m4^m\Gamma(\frac{2m+n+2}{2})\Gamma(\frac{2m-n+1}{2})}{\Gamma(\frac{n+2}{2})\Gamma(\frac{1-n}{2})}B_{0,n}(z), \\ &B_{2m+1,n}(z)=\frac{(-1)^m4^m\Gamma(\frac{2m+n+3}{2})\Gamma(\frac{2m-n+2}{2})}{\Gamma(\frac{n+3}{2})\Gamma(\frac{2-n}{2})}B_{1,n}(z)
\end{split}
\end{equation}
 and identical relations for $C_{m,n}(z)$. For $A_{2m,n}$ one gets:
\begin{equation}
\begin{split}
&A_{2m,n}(z)=\frac{(-1)^m4^m\Gamma(\frac{2m+n+2}{2})\Gamma(\frac{2m-n+1}{2})}{\Gamma(\frac{n+2}{2})\Gamma(\frac{1-n}{2})}A_{0,n}(z)+ \\ 
&+(-1)^{m+1}4^m\Gamma(\frac{2m+n+2}{2})\Gamma(\frac{2m-n+1}{2})\frac{z^{n+2}}{4}\sum_{j=0}^{m-1}(\frac{z}{2})^{2j}\frac{(-1)^j}{\Gamma(\frac{2k+n+4}{2})\Gamma(\frac{2k+3-n}{2})}
\end{split}
\end{equation}
whereas for $A_{2m+1,n}$ one has:
\begin{equation}
\begin{split}
&A_{2m+1,n}(z)=\frac{(-1)^m4^m\Gamma(\frac{2m+n+3}{2})\Gamma(\frac{2m-n+2}{2})}{\Gamma(\frac{n+3}{2})\Gamma(\frac{2-n}{2})}A_{1,n}(z)+ \\ 
&+(-1)^{m+1}4^m\Gamma(\frac{2m+n+3}{2})\Gamma(\frac{2m-n+2}{2})\frac{z^{n+3}}{4}\sum_{j=0}^{m-1}(\frac{z}{2})^{2j}\frac{(-1)^j}{\Gamma(\frac{2k+n+5}{2})\Gamma(\frac{2k+4-n}{2})}
\end{split}
\end{equation}

Let us keep separate the cases of $m$ even or odd.
\\
\textbf{The case $m$ even.} Notice that $n$ must be different from $m+2k+1$, $k=0,1,2,...$, so we must have $n<m$ (with $n \neq -m -2k, k=1,2,...$) or $n$ even. Let us assume $n$ even. Now we substitute $m \to 2m$ and $n \to 2n$ wherever in the previous results. We get 
\begin{equation}
\begin{split}
&B_{2m,2n}(z)=\frac{(-1)^m 4^{m}(m+n)!\Gamma(m-n+\frac{1}{2})}{n!\Gamma(\frac{1}{2}-n)}B_{0,2n}(z), \\ 
&C_{2m,2n}(z)=\frac{(-1)^m 4^{m}(m+n)!\Gamma(m-n+\frac{1}{2})}{n!\Gamma(\frac{1}{2}-n)}C_{0,2n}(z),
\end{split}
\end{equation}
$A_{0,2n}$, $B_{0,2n}$ and $C_{0,2n}$ are determined by the values of the functions $s_{\frac{1}{2},2n+\frac{1}{2}}(z)$. With reference to equation (\ref{pol}), we see that we must look at the following integral
\begin{equation}\label{mu0}
\int_0^1 \sin(zt)\frac{_2F_1\left(-2n,2n+1;1;\frac{1-t}{2}\right)}{_2F_1\left(-2n,2n+1;1;\frac{1}{2}\right)}dt.
\end{equation} 
The function $_2F_1\left(-2n,2n+1;1;\frac{1-t}{2}\right)$ is actually the Legendre polynomial $P_{2n}(t)$ of order $2n$ \cite{Pru}. Let us set
\begin{equation}
p_{2n}(t)\doteq \frac{P_{2n}(t)}{P_{2n}(0)}=\frac{_2F_1\left(-2n,2n+1;1;\frac{1-t}{2}\right)}{_2F_1\left(-2n,2n+1;1;\frac{1}{2}\right)}.
\end{equation}
Thanks to formula (\ref{intsin}) we get 
\begin{equation}\label{mu0}
z^{2n+1}\int \sin(zt)p_{2n}(t)dt=\sum_{j=0}^{n} (-1)^{j+1}z^{2n-2j}\cos(zt)p_{2n}^{2j}(t)+(-1)^j z^{2n-2j-1}\sin(zt)p_{2n}^{2j+1}(t)
\end{equation}
Evaluating the previous integral between $0$ and $1$ and confronting with equation (\ref{smn}) finally we get explicit expressions for $A_{0,2n}$, $B_{0,2n}$ and $C_{0,2n}$ in terms of $z$ with coefficients given in terms of derivatives of the Legendre polynomials $P_{2n}(t)$ of order $2n$ evaluated in $t=0$ or in $t=1$. More explicitly:
\begin{equation}\label{abc2n}
\begin{split}
& A_{0,2n}=\sum_{k=0}^n (-1)^k z^{2n-2k}p_{2n}^{2k}(0),\\
& B_{0,2n}=\sum_{k=0}^n (-1)^k z^{2n-2k}p_{2n}^{2k}(1),\\
& C_{0,2n}=\sum_{k=0}^{n-1} (-1)^{k+1} z^{2n-2k-1}p_{2n}^{2k+1}(1).
\end{split}
\end{equation}
In section (\ref{sec4}) a closed formula with explicit coefficients for these polynomials will be given. Now, we notice that the expression (\ref{smn}) for $m=0$ and $n \to 2n$, together with (\ref{mu0}) gives:
\begin{equation}
A_{0,2n}(z)-B_{0,2n}(z)\cos(z)-C_{0,2n}(z)\sin(z)=z^{2n+1}\int_0^{1} \sin(zt)p_{2n}(t)dt,
\end{equation} 
i.e.
\begin{equation}\label{order}
A_{0,2n}(z)-B_{0,2n}(z)\cos(z)-C_{0,2n}(z)\sin(z)=O(z^{2n+2}).
\end{equation}
In the theory of Pad\'e approximation \cite{Baker}, if two polynomials $q_n$ and $q_m$ of order, respectively, $n$ and $m$ are given and if they satisfy $q_n-f(z)q_m=O(z^{m+n+1})$, then  $q_m/q_n$ is a Pad\'e approximation of type $(m,n)$ of the function $f(z)$, i.e. the series of $q_m/q_n$ and the series for $f(z)$ agrees up to the order $m+n+1$. The equation (\ref{order})  can be interpreted in a similar manner: the polynomials $A_{0,2n}$, $B_{0,2n}$ and $C_{0,2n}$, of order respectively $2n$, $2n$ and $2n-1$ can be chosen in such a way that (\ref{order}) is satisfied. Notice that $A_{0,2n}$ and $B_{0,2n}$ are even polynomials whereas $C_{0,2n}(z)$ are odd polynomials and indeed their coefficients are not uniquely defined by (\ref{order}) as usually happens for Pad\'e approximants. By rewriting equation (\ref{order}) as
\begin{equation}\label{order1}
\frac{B_{0,2n}(z)}{A_{0,2n}(z)}\cos(z)+\frac{C_{0,2n}(z)}{A_{0,2n}(z)}\sin(z)=1+O(z^{2n+2}).
\end{equation}
it is clearer that $\frac{B_{0,2n}}{A_{0,2n}}$ and $\frac{C_{0,2n}}{A_{0,2n}}$ are Pad\'e approximants of the cosine and sine functions, i.e.
\begin{equation}\label{sc}
\frac{B_{0,2n}(z)}{A_{0,2n}(z)} = \cos(z)+O(z^{2n+2}), \quad \frac{C_{0,2n}(z)}{A_{0,2n}(z)}=\sin(z)+O(z^{2n+3}).
\end{equation}
From the previous it follows that the coefficients of the polynomials $A_{0,2n}$, $B_{0,2n}$ and $C_{0,2n}$ satisfy polynomials relations, since one has
\begin{equation}
B_{0,2n}(z)^2+C_{0,2n}(z)^2=A_{0,2n}(z)^2+O(z^{2n+2})
\end{equation}
The first few elements of the rational functions $\frac{B_{0,2n}}{A_{0,2n}} $ and $\frac{C_{0,2n}}{A_{0,2n}}$ are
\begin{equation}\begin{split}
&\frac{B_{0,2}}{A_{0,2}}=\frac{6-2z^2}{6+z^2}, \qquad \qquad\qquad \qquad\qquad \qquad\;\; \frac{C_{0,2}}{A_{0,2}}=\frac{6z}{6+z^2},\\
&\frac{B_{0,4}}{A_{0,4}}=\frac{840-360z^2+8z^4}{840+60z^2+3z^4}, \qquad\qquad\qquad \qquad  \frac{C_{0,4}}{A_{0,4}}=\frac{840z-80z^3}{840+60z^2+3z^4},\\
&\frac{B_{0,6}}{A_{0,6}}=\frac{166320-75600z^2+3360z^4-16z^6}{166320+7560z^2+210z^4+5z^6}, \quad  \frac{C_{0,6}}{A_{0,6}}=\frac{166320z-20160z^3+336z^5}{166320+7560z^2+210z^4+5z^6}.\\
\end{split}\end{equation}
\textbf{The case $m$ odd.} Again, $n$ must be different from $m+2k+1$, $k=0,1,2,...$, so we must have $n<m$ (with $n \neq -m -2k, k=1,2,...$) or $n$ odd. Let us assume $n$ odd. By substituting $m \to 2m+1$ and $n \to 2n+1$ we get 
\begin{equation}
\begin{split}
&B_{2m+1,2n+1}(z)=\frac{(-1)^m4^m(m+n+1)!\Gamma(m-n+\frac{1}{2})}{(n+1)!\Gamma(\frac{1}{2}-n)}B_{1,2n+1}(z), \\ 
&C_{2m+1,2n+1}(z)=\frac{(-1)^m4^m(m+n+1)!\Gamma(m-n+\frac{1}{2})}{(n+1)!\Gamma(\frac{1}{2}-n)}C_{1,2n+1}(z).
\end{split}
\end{equation}
$A_{1,2n+1}$, $B_{1,2n+1}$ and $C_{1,2n+1}$ are determined by the values of the functions $s_{\frac{3}{2},2n+\frac{1}{2}}(z)$. With reference to equation (\ref{pol}), we see that we must look at the following integral
\begin{equation}\label{mu0}
\int_0^1 (1-t)\sin(zt)\frac{_2F_1\left(-2n-1,2n+2;2;\frac{1-t}{2}\right)}{_2F_1\left(-2n-1,2n+2;2;\frac{1}{2}\right)}dt.
\end{equation} 
The polynomial $_2F_1\left(-2n-1,2n+2;2;\frac{1-t}{2}\right)$ is related to the associated Legendre polynomial $P_{2n+1}^{1}(t)$ of order 1 and degree $2n+1$ \cite{Pru}. Indeed one has
\begin{equation}
_2F_1\left(-2n-1,2n+2;2;\frac{1-t}{2}\right)=-\frac{1}{2(n+1)(2n+2)}\left(\frac{1+t}{1-t}\right)^{\frac{1}{2}}P_{2n+1}^{1}(t).
\end{equation}
More explicitly, these polynomials can be written in terms of derivatives of the function $(t^2-1)^{2n+1}$:
\begin{equation}
_2F_1\left(-2n-1,2n+2;2;\frac{1-t}{2}\right)=\frac{(1+t)}{2^{2n+2}(n+1)(2n+1)(2n+1)!}\frac{d^{2n+2}}{dt^{2n+1}}(t^2-1)^{2n+1}
\end{equation}
Let us define the polynomials $q_{2n+1}$ as:
\begin{equation}\label{qn}
q_{2n+1}(t)\doteq (1-t^2)^{\frac{1}{2}}\frac{P_{2n+1}^{1}(t)}{P_{2n+1}^{1}(0)}=\frac{(1-t) _2F_1\left(-2n-1,2n+2;2;\frac{1-t}{2}\right)}{_2F_1\left(-2n-1,2n+2;2;\frac{1}{2}\right)}.
\end{equation}
By noticing that $q_{2n+1}$ is a polynomial of degree $2n+2$, thanks to formula (\ref{intsin}) we get 
\begin{equation}\label{mu1}
z^{2n+3}\int \sin(zt)q_{2n+1}(t)dt=\sum_{j=0}^{n} (-1)^{j+1}z^{2n+2-2j}\cos(zt)q_{2n+1}^{2j}(t)+(-1)^j z^{2n+1-2j}\sin(zt)q_{2n}^{2j+1}(t),
\end{equation}
from which we obtain
\begin{equation}\label{a12n1}
\begin{split}
& A_{1,2n+1}=\sum_{k=0}^{n+1} (-1)^k z^{2n+2-2k}q_{2n+1}^{2k}(0),\\
& B_{1,2n+1}=\sum_{k=0}^{n+1} (-1)^k z^{2n+2-2k}q_{2n+1}^{2k}(1),\\
& C_{1,2n+1}=\sum_{k=0}^{n+1} (-1)^{k+1} z^{2n+1-2k}q_{2n+1}^{2k+1}(1).
\end{split}
\end{equation}
Also for these polynomials a closed formula will be given in section (\ref{sec4}). Again, one can look at the ratios $\frac{B_{1,2n+1}}{A_{1,2n+1}}$ and $\frac{C_{1,2n+1}}{A_{1,2n+1}}$ as Pad\'e approximant for the trigonometric functions
\begin{equation}\label{sc1}
\frac{B_{1,2n+1}(z)}{A_{1,2n+1}(z)} = \cos(z)+O(z^{2n+4}), \quad \frac{C_{1,2n+1}(z)}{A_{1,2n+1}(z)}=\sin(z)+O(z^{2n+5}).
\end{equation}
The first elements of the above rational functions are
\begin{equation}\begin{split}
&\frac{B_{1,1}}{A_{1,1}}=\frac{2}{2+z^2}, \qquad \qquad\qquad \qquad\quad\;\;\; \frac{C_{1,1}}{A_{1,1}}=\frac{2z}{2+z^2},\\
&\frac{B_{1,3}}{A_{1,3}}=\frac{120-48z^2}{120+12z^2+z^4}, \qquad\qquad \quad \;\; \frac{C_{1,3}}{A_{1,3}}=\frac{120z-8z^3}{120+12z^2+z^4},\\
&\frac{B_{1,5}}{A_{1,5}}=\frac{15120-6720z^2+240z^4}{15120+840z^2+30z^4+z^6}, \quad  \frac{C_{1,5}}{A_{1,5}}=\frac{15120z-1680z^3+16z^5}{15120+840z^2+30z^4+z^6}.\\
\end{split}\end{equation}

In general, the polynomials defined in (\ref{smn}) are Pad\'e approximant for the trigonometric functions. Indeed, from (\ref{smn}) and (\ref{intsin}) it follows 
\begin{equation}\label{sc1}
\frac{B_{m,n}(z)}{A_{m,n}(z)} = \cos(z)+O(z^{m+n+2}), \quad \frac{C_{m,n}(z)}{A_{m,n}(z)}=\sin(z)+O(z^{m+n+3}).
\end{equation}

\section{Explicit formulae and some numerics}\label{sec4}
\begin{figure}[H]
\centering
\includegraphics[scale=0.5]{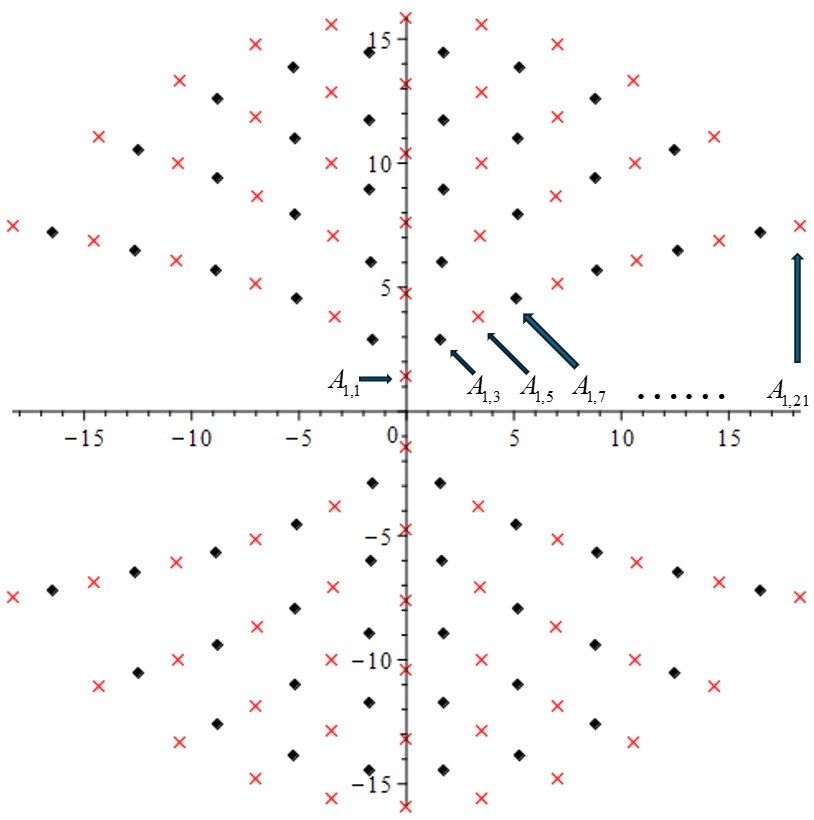}
\caption{The zeros of the polynomials $A_{1,2n+1}$ for $n= 1,2,...,10$.}
\label{fig1}
\end{figure}

\begin{figure}[H]
\centering
\includegraphics[scale=0.5]{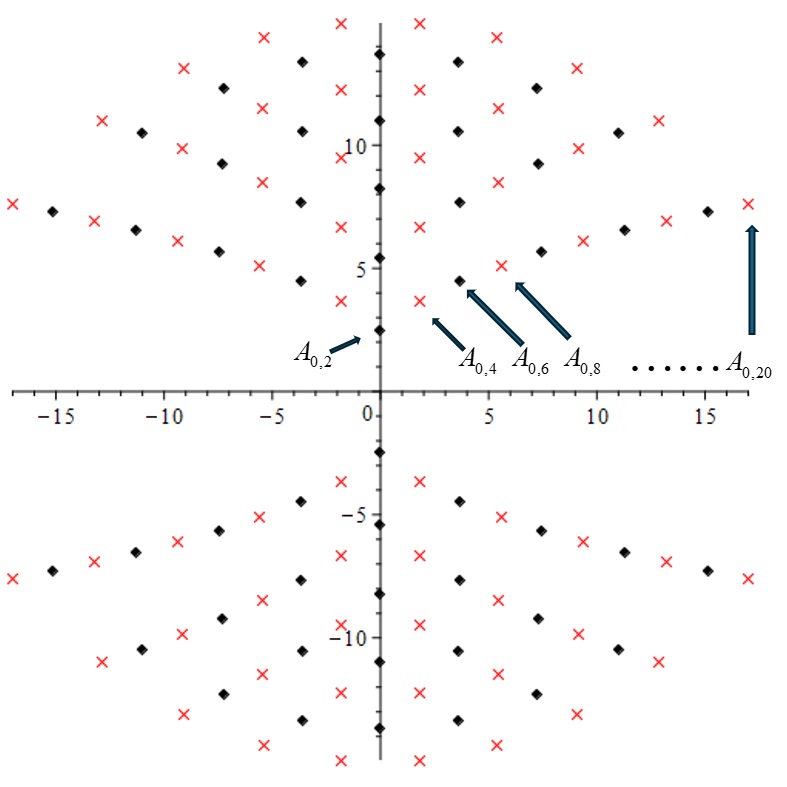}
\caption{The zeros of the polynomials $A_{0,2n}$ for $n= 1,2,...,10$.}
\label{fig2}
\end{figure}

The rational functions $\frac{B_{0,2n}(z)}{A_{0,2n}(z)} $, $\frac{C_{0,2n}(z)}{A_{0,2n}(z)} $, $\frac{B_{1,2n+1}(z)}{A_{1,2n+1}(z)} $ and $\frac{C_{1,2n+1}(z)}{A_{1,2n+1}(z)} $ are Pad\'e approximants for the trigonometric functions $\cos(z)$ and $\sin(z)$. It is natural to ask about how distribute the zeros of the numerators and denominators of these polynomials. We expect that the zeros of $B_{m,n}$ and $C_{m,n}$ are close to the zeros of $\cos(z)$ and $\sin(z)$, whereas $A_{m,n}$ is expected to have no real zeros. Indeed, we can give a closed formula for the coefficients of these polynomials. Let us start with $A_{0,2n}$. Since the Legendre polynomials of degree $2n$ are explicitly given by
\begin{equation}
P_{2n}(t)=\frac{1}{2^{2n}}\sum_{k=0}^n \frac{(-1)^k(4n-2k)!}{k!(2n-k)!(2n-2k)!}t^{2n-2k},
\end{equation}
by using formula (\ref{abc2n}), making the corresponding derivatives, evaluating them to $0$, after some manipulations we get:
\begin{equation}\label{a02n}
A_{0,2n}(z)=\frac{(n!)^2}{(2n)!}\sum_{k=0}^n \frac{(2n+2k)!}{(n+k)!(n-k)!}z^{2n-2k}.
\end{equation}
From the above formula we see that all the coefficients are positive. Also, the polynomials $A_{0,2n}$ are even in $z$: they have no real roots. For the functions $B_{0,2n}$ and $C_{0,2n}$ we can use the expansion of $P_{2n}(t)$ around $t=1$, i.e.
\begin{equation}\label{pin1}
P_{2n}(t)=\sum_{k=0}^{2n} \frac{(-1)^k(2n+k)!}{(k!)^2(2n-k)!}\left(\frac{1-t}{2}\right)^k.
\end{equation} 
By looking at (\ref{abc2n}) and making the corresponding derivatives of (\ref{pin1}) we get:
\begin{equation}\label{b02n}
B_{0,2n}(z)=\frac{(n!)^2(-1)^n}{(2n)!}\sum_{k=0}^n (-1)^k\frac{(2n+2k)!}{(2k)!(2n-2k)!}(2z)^{2n-2k}.
\end{equation}
Analogously, for $C_{0,2n}$, we get
\begin{equation}\label{c02n}
C_{0,2n}(z)=\frac{(n!)^2(-1)^n}{(2n)!}\sum_{k=0}^{n-1} (-1)^{k+1}\frac{(2n+2k+1)!}{(2k+1)!(2n-2k-1)!}(2z)^{2n-2k-1}.
\end{equation} 
For the polynomials $A_{1,2n+1}$, $B_{1,2n+1}$ and $C_{1,2n+1}$ we could use equations (\ref{qn}) and (\ref{a12n1}) or, more easily, the recurrences (\ref{recu}), from which it follows that 
\begin{equation}
A_{2n+1}=\frac{(2n+1)A_{0,2n}+2z^2(n+1)A_{0,2n+2}}{4n+3},
\end{equation}
and equal identities for $C_{1,2n+1}$ and $B_{1,2n+1}$. By using (\ref{a02n}), (\ref{b02n}) and (\ref{c02n}) we get:
\begin{equation}
A_{1,2n+1}=\frac{4(2n+1)((n+1)!)^2}{(2n+2)!}\sum_{k=0}^{n+1} \frac{(2n+2k-1)!}{(n+k-1)!(n-k-1)!}z^{2n-2k+2},
\end{equation}
\begin{equation}
B_{1,2n+1}(z)=\frac{2(2n+1)((n+1)!)^2(-1)^n}{(2n+2)!}\sum_{k=0}^n (-1)^k\frac{(2n+2k+2)!}{(2k+1)!(2n-2k)!}(2z)^{2n-2k},
\end{equation}
\begin{equation}
C_{1,2n+1}(z)=\frac{2(2n+1)((n+1)!)^2(-1)^n}{(2n+2)!}\sum_{k=0}^n (-1)^k\frac{(2n+2k+1)!}{(2k)!(2n-2k+1)!}(2z)^{2n-2k+1}.
\end{equation}
Numerically, we have seen that, for $n\in (0,10)$, indeed $B_{0,2n}$, $C_{0,2n}$, $B_{1,2n+1}$ and $C_{1,2n+1}$ possess only real zeros and, for $z$ relatively small, the zeros approximate very well the values of the zeros of $\cos(z)$ and $\sin(z)$. In table (\ref{t1}) we report the relative difference  $\frac{z_{k}^{n}-n\pi}{n\pi}$ between the $n^{th}$ positive zero $z_{k}^n$ of the polynomial $C_{0, 2k}(z)$ and the $n^{th}$ positive zero of the $\sin$ function $n\pi$ for $k=1..6$.  In table (\ref{t2}) the same for the zeros of $B_{1,2k+1}(z)$ compared to the zeros of $\cos(z)$. Also, in figures (\ref{fig1}) and (\ref{fig2}) we report respectively the zeros of $A_{0,2n}$ and $A_{1,2n+1}$ for $n=0..10$: these appears to form a pattern in the complex plane.

\begin{table}[H]
\begin{center}
\caption{Values of the relative difference $\frac{z_{k}^{n}-n\pi}{n\pi}$ between the $n^{th}$ positive zero $z_{k}^n$ of the polynomial $C_{0, 2k}(z)$ and the $n^{th}$ positive zero of the $\sin$ function $n\pi$.}
\label{t1}
\begin{tabular}{|c|c|c|c|c|c|c|c|}
\hline
 \diagbox{$k$}{$n$} & 1 & 2 & 3 & 4 & 5 & 6 \\ \hline
 1 & $3.14\cdot 10^{-2}$ & - & - & - & - & -    \\ \hline
  2 & $2.93\cdot 10^{-4}$  & $1.27\cdot 10^{-1}$ & - & - & - & -  \\ \hline
 3 & $7.07\cdot 10^{-7}$ & $6.56\cdot 10^{-3}$ & $2.67\cdot 10^{-1}$ & - & - & -   \\ \hline
 4 & $6.03\cdot 10^{-10}$ & $1.37\cdot 10^{-4}$ & $2.94\cdot 10^{-2}$ & $4.33\cdot 10^{-1}$ & - & -  \\ \hline
  5 & $2.29\cdot 10^{-13}$ & $1.10\cdot 10^{-6}$ & $1.90\cdot 10^{-3}$ & $7.26\cdot 10^{-2}$ & $6.13\cdot 10^{-1}$ & -    \\ \hline
 6 & $4.46\cdot 10^{-17}$ & $4.15\cdot 10^{-9}$ & $5.29\cdot 10^{-5}$ & $9.27\cdot 10^{-3}$ & $1.33\cdot 10^{-1}$ &  $8.01\cdot 10^{-1}$    \\ \hline 
\end{tabular}
\end{center}
\end{table}

\begin{table}[H]
\begin{center}
\caption{Values of the relative difference $\frac{z_{k}^{n}-\left(\frac{\pi}{2}+(n-1)\pi\right)}{\frac{\pi}{2}+(n-1)\pi}$ between the $n^{th}$ positive zero $z_{k}^n$ of the polynomial $B_{1, 2k+1}(z)$ and the $n^{th}$ positive zero of the $\cos$ function $\frac{\pi}{2}+(n-1)\pi$.}
\label{t2}
\begin{tabular}{|c|c|c|c|c|c|c|c|}
\hline
 \diagbox{$k$}{$n$} & 1 & 2 & 3 & 4 & 5 & 6 \\ \hline
 1 & $6.58\cdot 10^{-3}$ & - & - & - & - & -    \\ \hline
  2 & $7.36\cdot 10^{-6}$  & $7.23\cdot 10^{-2}$ & - & - & - & -  \\ \hline
 3 & $1.72\cdot 10^{-9}$ & $1.95\cdot 10^{-4}$ & $1.93\cdot 10^{-1}$ & - & - & -   \\ \hline
 4 & $1.32\cdot 10^{-13}$ & $1.70\cdot 10^{-5}$ & $1.54\cdot 10^{-2}$ & $3.48\cdot 10^{-1}$ & - & -  \\ \hline
  5 & $4.26\cdot 10^{-18}$ & $5.50\cdot 10^{-8}$ & $6.13\cdot 10^{-4}$ & $4.85\cdot 10^{-2}$ & $5.22\cdot 10^{-1}$ & -    \\ \hline
 6 & $6.79\cdot 10^{-23}$ & $8.16\cdot 10^{-11}$ & $9.89\cdot 10^{-6}$ & $4.59\cdot 10^{-3}$ & $0.11\cdot 10^{-1}$ &  $0.71\cdot 10^{-1}$    \\ \hline 
\end{tabular}
\end{center}
\end{table}

\section{Trigonometric representation for $\mu$ an integer.}
When $\mu$ is a positive integer it is possible to give explicit formulae for the hypergeometric kernel in terms of trigonometric functions. We start from the formula (\ref{int}) that can be rewritten also in terms of an integral involving the $\sin(z\cos(\theta))$ in the following way:
\begin{equation}\label{intnew}
s_{0,\nu}(z)=\frac{1}{\cos\left(\frac{\nu\pi}{2}\right)}\int_0^\frac{\pi}{2} \sin(z\cos(\theta))\cos(\nu \theta)d\theta.
\end{equation}
Since the functions $s_{n+1,\nu}(z)$ solve the recurrence (\ref{recu})
\begin{equation}\label{rec}
\frac{2\nu}{z}s_{n+1,\nu}(z)=(n+\nu)s_{n,\nu-1}(z)-(n-\nu)s_{n,\nu+1}(z),
\end{equation}
immediately one has the corresponding integral representation for $s_{1,\nu}(z)$: 
\begin{equation}
s_{1,\nu}(z)=\frac{z}{\sin\left(\frac{\nu\pi}{2}\right)}\int_{0}^{\frac{\pi}{2}}\sin(z\cos(\theta))\sin(\nu \theta)\sin(\theta)d\theta.%=\frac{z}{\sin\left(\frac{\nu\pi}{2}\right)}\int_0^1\sin(zt)\sin(\nu\arccos(t))dt
\end{equation}
%The previous equation can be checked directly by using the series of $\sin(z\cos(t))$ and the equivalence (see e.g. \cite{Ryzhik}, equation 3.633.1)
%\begin{equation}
%\int_{0}^{\frac{\pi}{2}}\cos(t)^{2k+1}\sin(\nu t)\sin(t)dt=\frac{\nu\pi}{2^{2k+4}(k+1)}\frac{(2k+2)!}{\Gamma(k+2+\frac{\nu}{2})\Gamma(k+2-\frac{\nu}{2})},
%\end{equation}
%giving
%\begin{equation}
%s_{1,\nu}(z)=\sum_{k=0}c_kz^{2k+2},
%\end{equation}
%with
%\begin{equation}
%c_k=\frac{(-1)^k\nu\pi}{2^{2k+3}\sin\left(\frac{\nu\pi}{2}\right)}\frac{1}{\Gamma(k+2+\frac{\nu}{2})\Gamma(k+2-\frac{\nu}{2})}=\frac{(-1)^k}{\prod_{m=1}^{k+1}(2m+\nu)(2m-\nu)}.
%\end{equation}
It is clear from (\ref{rec}) that, in general, we can represent the entire function $s_{n,\nu}(z)$ as an integral of the type:
 \begin{equation}\label{ans1}
s_{n,\nu}(z)=z^n\int_{0}^{\frac{\pi}{2}}\sin(z\cos(\theta))f_n(\nu,\theta)\sin(\theta)d\theta,
\end{equation}
where $f_n(\nu,\theta)$ is a linear combination of $\sin$ functions. More explicitly, from (\ref{rec}) we get:
\begin{equation}\label{ans2}
f_n(\nu,\theta)=\sum_{k=0}^{n-1}a_{k}^n(\nu)\frac{\sin\left((\nu-n+2k+1)\theta\right)}{\sin\left((\nu-n+2k+1)\frac{\pi}{2}\right)}, \quad n\geq 1.
\end{equation}
The coefficients $a_{k}^n(\nu)$ are rational functions of $\nu$ and solve a recurrence in $n$ and $\nu$ that follows from (\ref{rec}). However, to get an explicit expression for these coefficients, we can directly insert the ansatze (\ref{ans1}) and (\ref{ans2}) in (\ref{Lom}) to obtain a different recursion for $a_k^n(\nu)$. Indeed, the functions (\ref{ans1}) are solutions of (\ref{Lom}) provided that the functions $f_n(\nu,\theta)$ solve the following differential equation:
\begin{equation}\label{feq}
\sin(\theta)\frac{d^2f_n(\nu,\theta)}{d\theta^2}-2(n-1)\cos(\theta)\frac{df_n(\nu,\theta)}{d\theta}-\sin(\theta)\left((n-1)^2-\nu^2\right)f_n(\nu,\theta)=0,
\end{equation}
with the boundary condition $f_n(\nu,\frac{\pi}{2})=1$ (the other necessary boundary condition $f_n(\nu,0)=0$ is automatically satisfied by the ansatz (\ref{ans2})). By inserting (\ref{ans1}) and (\ref{ans2}) in (\ref{feq}), after few manipulation we get that equation (\ref{feq}) is satisfied if the following recursion for $a_k^n(\nu)$ holds:
\begin{equation}
k(k+\nu)a_k^n(\nu)+(k-n)(k-n+\nu)a_{k-1}^n(\nu)=0, \quad k=1..n-1.
\end{equation} 
This recurrence gives the value of $a_k^{n}(\nu)$, $k=1..n-1$, in terms of $a_0^n(\nu)$. The value of  $a_0^n(\nu)$ then must be fixed to satisfy the boundary condition $f_n(\nu,\frac{\pi}{2})=1$. We get
\begin{equation}
a_{0}^n(\nu)=\frac{1}{2^{n-1}}\prod_{p=1}^{n-1}\frac{\nu-n+2p+1}{\nu-p+1},
\end{equation}
and the following explicit formula for $a_k^n(\nu)$:
\begin{equation}\label{akn}
a_k^n(\nu)=\frac{1}{2^{n-1}}\prod_{p=1}^{n-1}\frac{\nu-n+2p+1}{\nu-p+1}\prod_{q=1}^k\frac{(n-q)(q+\nu-n)}{q(q+\nu)}.
\end{equation}

We notice that from formulae (\ref{ans1}), (\ref{ans2}) and (\ref{akn}), by comparing with formula (\ref{intr}) and with some manipulations we get the following expression relating the hypergeometric functions $ _2F_1$ with trigonometric functions for $n\geq 1$ and $|\theta|<\pi$:
\footnotesize{\begin{equation}\label{trigexp}
_2F_1\left(\frac{1}{2}+\nu,\frac{1}{2}-\nu;n+\frac{1}{2};\sin(\frac{\theta}{2})^2\right)=\frac{\prod_{p=0}^{n-1}\frac{2p+1}{\nu-p}}{2^{3n-2}\sin(\frac{\theta}{2})^{2n-1}}\sum_{k=0}^{n-1}\prod_{q=1}^k\frac{(q-n)(q+\nu-n)}{q(q+\nu)}\sin\left((1-n+\nu+2k)\theta\right) 
\end{equation}}
\normalsize
The equivalence (\ref{trigexp}) generalizes the equivalent expression for $n=1$ given in \cite{GR} (see equation 9.121-16), see also \cite{Chu}. The first few examples are the following:
\footnotesize{\begin{equation}
\begin{split}
& _2F_1\left(\frac{1}{2}+\nu,\frac{1}{2}-\nu;\frac{3}{2};\sin(\frac{\theta}{2})^2\right)=\frac{\sin(\nu\theta)}{2\nu\sin(\frac{\theta}{2})},\\
& _2F_1\left(\frac{1}{2}+\nu,\frac{1}{2}-\nu;\frac{5}{2};\sin(\frac{\theta}{2})^2\right)=\frac{3}{16\nu\sin(\frac{\theta}{2})^3}\left(\frac{\sin((\nu-1)\theta)}{(\nu-1)}-\frac{\sin((\nu+1)\theta)}{(\nu+1)}\right),\\
& _2F_1\left(\frac{1}{2}+\nu,\frac{1}{2}-\nu;\frac{7}{2};\sin(\frac{\theta}{2})^2\right)=\frac{15}{128\nu\sin(\frac{\theta}{2})^5}\left(\frac{\sin((\nu-2)\theta)}{(\nu-2)(\nu-1)}-\frac{2\sin(\nu\theta)}{(\nu-1)(\nu+1)}+\frac{\sin((\nu+2)\theta)}{(\nu+1)(\nu+2)}\right),\\
&\cdots
\end{split}
\end{equation}}
\normalsize
Finally, let us take $\nu=\mu$ in (\ref{intr}): we get the well known representation for $s_{\mu,\mu}(z)$  ($s_{\mu,\mu}(z)$ is proportional to the so-called Struve function):
 \begin{equation}\label{intr1}
s_{\mu,\mu}(z)=z^{\mu}\int_0^1 (1-t^2)^{\mu-\frac{1}{2}}\sin(zt)dt,
\end{equation}
We notice also that when $\nu=\mu+2n$, with $n\in \mathbb{N}$, it is possible to explicitly write the integral (\ref{intr}) in terms of an algebraic kernel, like in (\ref{intr0}) and (\ref{intr1}). The first two terms are
 \begin{equation}\label{intr2}
s_{\mu,\mu+2}(z)=z^{\mu}\int_0^1 (1-t^2)^{\mu-\frac{1}{2}}\left(1-2(\mu+1)t^2\right)\sin(zt)dt,
\end{equation}
and
 \begin{equation}\label{intr3}
s_{\mu,\mu+4}(z)=z^{\mu}\int_0^1 (1-t^2)^{\mu-\frac{1}{2}}\left(1-4(\mu+2)t^2+\frac{4}{3}(\mu+2)(\mu+3)t^4\right)\sin(zt)dt.
\end{equation}
\\
\begin{center} {\bf Acknowledgments} \end{center}
I wish to acknowledge the support of Universit\`a degli Studi di Brescia,  GNFM-INdAM and INFN, Gr. IV - Mathematical Methods in NonLinear Physics. Further, I would like to acknowledge the support of ISNMP - International Society of Nonlinear Mathematical Physics.


\normalsize





\small 
\begin{thebibliography}{40}
\bibitem{Baker} Baker, G. A. Jr., Graves-Morris, P.: Pad\'e Approximants. New York: Cambridge University Press, 1996
\bibitem{Chu} W. Chu: Trigonometric expressions for Gaussian $_2F_1$ series, Turk. J. Math., 43 (4) (2019), pp. 1823-1836
\bibitem{GR} I.S. Gradshteyn, I.M. Ryzhik: Table of integrals, series, and products, 7th edition, D. Zwillinger and A. Jeffrey Ed, Academic Press, Amsterdam, 2007.
\bibitem{KL} Koumandos, S., Lamprecht, M.: The zeros of certain Lommel functions. Proc.
Am. Math. Soc. 140(9), 3091-3100 (2012)

\bibitem{NIST} NIST Digital Library of Mathematical Functions. https://dlmf.nist.gov/, Release 1.1.11 of
2023-09-15. F. W. J. Olver, A. B. Olde Daalhuis, D. W. Lozier, B. I. Schneider, R. F. Boisvert,
C. W. Clark, B. R. Miller, B. V. Saunders, H. S. Cohl, and M. A. McClain, eds.
\bibitem{P} G. P\'olya: \"Uber die Nullstellen gewisser ganzer Funktionen, Mathematische Zeitschrift, 2, 352–383, 1918.
\bibitem{Pru} A. P. Prudnikov, Yu. A. Brychkov, and O. I. Marichev: Integrals and Series: More Special Functions, Vol. 3. Translated by G.G. Gould. Gordon and Breach Science Publishers, New York, 1990.


\bibitem{W} G. N. Watson: A Treatise on the Theory of Bessel Functions, Cambridge University Press, London, 1922.
\bibitem{Z} F. Zullo: Notes on the zeros of the solutions of the non-homogeneous Airy's equation, in  Formal and Analytic Solutions of Differential Equations, G. Filipuk, A. Lastra and S. Michalik eds., pp. 125-144, World Scientific Publishing Europe Ltd., 2022.
\bibitem{Z1} F. Zullo: Integral representations and zeros of the Lommel function and the hypergeometric $ \,_1F_2 $ function. Submitted

\end{thebibliography}
\end{document}